\documentclass[11pt,oneside,a4paper]{amsart}
\usepackage[T1]{fontenc}
\usepackage[utf8]{inputenc}
\usepackage[margin=1.5in]{geometry}
\usepackage{latexsym,amssymb,upref,amsmath,amssymb,amsfonts,amsthm}
\usepackage{graphicx}
\usepackage[nobysame,msc-links]{amsrefs} 
\usepackage{upref}
\usepackage[english]{babel}
\usepackage{hyperref} 

\newcommand{\articletitle}{Terminal-Pairability in Complete Bipartite Graphs}

\hypersetup 
{
  unicode=true,          
  pdfstartview={XYZ null null 1},    
  pdftitle={\articletitle},
  pdfauthor={Lucas Colucci, P\'eter L. Erd\H{o}s, Ervin Gy\H{o}ri, Tam\'as R\'obert Mezei},
  pdfsubject={Mathematics, Graph theory},   
  pdfkeywords={Terminal-pairability, Complete bipartite graphs}, 
	colorlinks=true,        
  citecolor=red,linkcolor=red,urlcolor=blue
}

\newtheorem{theorem}{Theorem}
\newtheorem{conj}[theorem]{Conjecture}

\newcommand\blfootnote[1]{%
  \begingroup
  \renewcommand\thefootnote{}\footnote{#1}%
  \addtocounter{footnote}{-1}%
  \endgroup
}

\makeatletter
\g@addto@macro{\endabstract}{\@setabstract}
\newcommand{\authorfootnotes}{\renewcommand\thefootnote{\@fnsymbol\c@footnote}}%
\makeatother
\begin{document}
\thispagestyle{empty}
\addtocounter{footnote}{2}

\mbox{}
\begin{center}
  \vspace{36pt}

  \Large
  {\bfseries \articletitle} \par \bigskip

  \blfootnote{E-mail addresses: \texttt{lucas.colucci.souza@gmail.com}, \texttt{erdos.peter@renyi.mta.hu},\\\texttt{gyori.ervin@renyi.mta.hu}, \texttt{tamasrobert.mezei@gmail.com}}
  \authorfootnotes
  \small\scshape
  Lucas Colucci\textsuperscript{2},
  P\' eter L.~Erd\H os\textsuperscript{1,}\footnote{Research  of the author was supported by the National Research, Development and Innovation --- NKFIH grant K~116769.}\textsuperscript{,}\footnote{Research  of the author was supported by the National Research, Development and Innovation --- NKFIH grant SNN~116095.}, \\
  Ervin Gy\H ori\textsuperscript{1,2,}\addtocounter{footnote}{-2}\footnotemark\textsuperscript{,}\footnotemark,
  Tam\'as R\'obert Mezei\textsuperscript{1,2,}\addtocounter{footnote}{-2}\footnotemark\textsuperscript{,}\addtocounter{footnote}{1}\footnote{Corresponding author}
  \par \bigskip

  \normalfont
  \textsuperscript{1}Alfréd Rényi Institute of Mathematics, Hungarian Academy of Sciences, Re\'altanoda~u.~13--15, 1053 Budapest, Hungary \par
  \textsuperscript{2}Central European University, Department of Mathematics and its Applications, N\'ador~u.~9, 1051 Budapest, Hungary \par \bigskip

  \normalsize\today
\end{center}

\begin{abstract}
  We investigate the terminal-pairibility problem in the case when the base graph is a complete bipartite graph, and the demand graph is also bipartite with the same color classes. We improve the lower bound on maximum value of $\Delta(D)$ which still guarantees that the demand graph $D$ is terminal-pairable in this setting. We also prove a sharp theorem on the maximum number of edges such a demand graph can have.
\end{abstract}

\section{Introduction}

The \emph{terminal-pairability} problem has been introduced in~\cite{TPP92}. It asks the following question: given a simple \emph{base graph} $G$ and a list of pairs of vertices of $G$ (which list may contain multiple copies of the same pair), can we assign to each pair a path in $G$ whose end-vertices are the two elements of the pair, such that the set of chosen paths are pairwise edge-disjoint.

\medskip

The above problem can be compactly described by a pair of graphs: the base graph and a so-called \emph{demand graph}, which is a loopless multigraph on the same set of vertices as the base graph together with the list of pairs to be joined as the (multi)set of edges. If the terminal-pairability problem defined by $D$ and $G$ can be solved, then we say that $D$ is resolvable in $G$. In this paper, demand graphs are denoted by $D$, or its primed and/or indexed variants.

\medskip

Related to the terminal-pairability problem is the notion of weak linkedness, which is closely tied to the edge-connectivity number (see~\cites{MR1143181}). A graph $G$ is weakly-$k$-linked if and only if every demand graph on $V(G)$ with at most $k$ edges is resolvable in $G$. In the terminal-pairability context, however, we are more interested in the degrees of $D$.

\medskip

Given an edge $e\in E(D)$ with endvertices $x$ and $y$, we define the \emph{lifting of $e$ to a vertex $z\in V(D)$}, as an operation which transforms $D$ by deleting $e$ and adding two new edges joining $xz$ and $zy$; in case $z=x$ or $z=y$, the operation does not do anything. We stress that we do not use any information about $G$ to perform a lifting and that the graph obtained using a lifting operation is still a demand graph.

\medskip

Notice that the terminal-pairability problem defined by $G$ and $D$ is solvable if and only if there exists a series of liftings, which, applied successively to $D$, results in a (simple!) subgraph of $G$. This subgraph is called a \emph{resolution} of $D$ in $G$. The edge-disjoint paths can be recovered by assigning pairwise different labels to the edges of $D$, and performing the series of liftings so that new edges inherit the label of the edge they replace. Clearly, edges sharing the same label form a walk between the endpoints of the demand edge of the same label in $D$, and so there is also such a path.

\medskip

This problem has been studied, for example, for complete graphs \cites{TPP92, GyMM16Complete} and cartesian product of complete graphs~\cites{KKGL99, GyMM16Grid}. In this paper we deal with problems where the base graph is a complete bipartite graph and the demand graph is bipartite with the color classes of the base graph.

\begin{conj}[\cite{GyarfasSchelp98}]\label{conj:nover3}
  Let $D$ be a bipartite demand graph whose base graph is $K_{n,n}$, i.e., $V(D)=V(K_{n,n})$ and each element of $E(D)$ is a copy of an edge of $K_{n,n}$. If $\Delta(D)\le \lceil n/3\rceil $ holds, then $D$ is resolvable in $K_{n,n}$.
\end{conj}

The above conjecture is sharp in the sense that the disjoint union of $n$ pairs of vertices each joined by $\lceil n/3\rceil+1$ parallel edges cannot be resolved in $K_{n,n}$, as explained by the following reasoning. From each set of edges joining the same pair of vertices at most one edge is resolved into a path of length 1 (itself), while the rest of them must be replaced by paths of length at least 3, therefore any resolution uses at least $n+3\cdot n\cdot  \lceil n/3\rceil \ge n^2+n$ edges in $K_{n,n}$, which is a contradiction.

\medskip

By replacing $\lceil n/3\rceil$ with $n/12$ in Conjecture~\ref{conj:nover3}, we get a theorem of Gyárfás and Schelp~\cite{GyarfasSchelp98}. We also cannot prove Conjecture~\ref{conj:nover3} in its generality, but in the following theorem we improve the previous best known bound of $n/12$ to $(1-o(1))n/4$.

\begin{theorem}\label{thm:deg}
  Let $D$ be a bipartite demand graph whose two color classes $A$ and $B$ have sizes $a$ and $b$, respectively.
  If $d(x)\le (1-o(1))b/4$ for all $x\in A$ and $d(y)\le (1-o(1))a/4$ for all $y\in B$, then $D$ is resolvable in the complete bipartite graph with color classes $A$ and $B$.
\end{theorem}

For certain graph classes, if $n$ is divisible by 3, we can prove that the sharp bound $n/3$ holds. Let $\uplus$ denote the disjoint union of sets.

\begin{theorem}\label{thm:optcase}
  Let $D$ be a bipartite demand graph with base graph $K_{n,n}$, such that
  \[U=\biguplus\limits_{i=1}^3 U_i\text{ and }V=\biguplus\limits_{i=1}^3 V_i\] are the two color classes of $D$ with $|U_i|=|V_i|\ge \lfloor\frac{n}{3}\rfloor$ for $i=1,2,3$.
  If $\Delta(D)\le \lfloor\frac{n}{3}\rfloor$ and for any $i\neq j$ there is no edge of $D$ joining some vertex of $U_i$ to some vertex of $V_j$,
  then $D$ is resolvable in $K_{n,n}$.
\end{theorem}

Additionally, we prove a sharp bound on the maximum number of edges in a resolvable bipartite demand graph:

\begin{theorem}\label{thm:edge}
  Let $n\geq 4$ and $D$ be a bipartite demand graph with the base graph $K_{n,n}$. If $D$ has at most $2n-2$ edges and $\Delta(D)\le n$, then $D$ is resolvable in $K_{n,n}$.
\end{theorem}

Notice the assumption $\Delta(D)\le n$ is necessary: there can be at most $n$ edge-disjoint paths starting at any given vertex. The result is sharp, as it is shown by the demand graph composed of a pair of vertices joined by $n$ edges, another pair of vertices joined by $n-1$ edges, and $2n-4$ isolated vertices: in any resolution, one of the paths corresponding to one of the $n$ edges joining the first pair of vertices passes through a vertex of the pair of vertices joined by $n-1$ edges, implying that this vertex has degree $\ge n+1$ in the resolution, a contradiction.

\section{Proofs of the degree versions (Theorem~\ref{thm:deg} and~\ref{thm:optcase})}

Theorem~\ref{thm:optcase} serves a dual purpose in our analysis: it provides several examples where Conjecture~\ref{conj:nover3} holds and it demonstrates the techniques that will be used in the proof of Theorem~\ref{thm:deg}.
Before we proceed to prove the theorems, we state several definitions and three well-known results about edge-colorings of multigraphs.

\medskip

Let $H$ be a loopless multigraph. Recall that the chromatic index (or the edge chromatic number) $\chi'(H)$ is the minimum number of colors required to properly color the edges of a graph $H$. Similarly, the list chromatic index (or the list edge chromatic number) $\mathrm{ch}'(H)$ is the smallest integer $k$ such that if for each edge of $G$ there is a list of $k$ different colors given, then there exists a proper coloring of the edges of $H$ where each edge gets its color from its list. The maximum multiplicity $\mu(H)$ is the maximum number of edges joining the same pair of vertices in $H$. The number of edges joining a vertex $x\in V(H)$ to a subset $A\subseteq V(H)$ of vertices is denoted by $e_H(x,A)$. The set of neighbors of $x$ in $H$ is denoted by $N_H(x)$. For other notation the reader is referred to~\cite{Diestel}.

\medskip

\begin{theorem}[Kőnig~\cite{Konig16}]
  For any bipartite multigraph $H$ we have $\chi'(H)=\Delta(H)$, or, in other words, the edge set of $H$ can be decomposed into $\Delta(H)$ matchings.
\end{theorem}

\medskip

\begin{theorem}[Vizing,~\cite{Vizing65}]\label{thm:vizing}
  For any multigraph $H$
  \[\chi'(H)\le \Delta(H)+\mu(H).\]
\end{theorem}

\medskip

\begin{theorem}[Kahn,~\cite{Kahn00}]\label{thm:kahn}
  For any multigraph $H$
  \[\mathrm{ch}'(H)\le (1+o(1))\chi'(H).\]
\end{theorem}

\medskip

Even though in our theorems the demand graphs are bipartite, in the proofs we may transform them into non-bipartite ones.

\begin{proof}[\normalfont\bfseries Proof of Theorem~\ref{thm:optcase}]
  Let $D_i$ be the (bipartite) subgraph of $D$ induced by $U_i\cup V_i$ for $i=1,2,3$. As parallel edges are allowed in $D$, without loss of generality, we may assume that $D_i$ is $\lfloor\frac{n}{3}\rfloor$-regular.
  By Kőnig's theorem, $E(D_i)$ can be partitioned into matchings $M_{i,1},M_{i,2},\ldots,M_{i,\lfloor\frac{n}{3}\rfloor}$, each of size $|U_i|$.
  We derive $D'_i$ from $D$ by lifting the edges of $M_{i,j}$ to the $j$\textsuperscript{th} vertex of $U_i$ for each $j=1,2,\ldots,\lfloor\frac{n}{3}\rfloor$.
  Firstly, all the edges of $D'_i$ between $U_i$ and $V_i$ have multiplicity 1.
  Secondly, observe, that $D'_i[U_i]$ is $2(\lfloor\frac{n}{3}\rfloor-1)$-regular and $\mu(D'_i[U_i])=2$.

  \medskip

  Applying Vizing's theorem we get $\chi'(D'_i[U_i])\le \Delta(D'_i[U_i])+\mu(D'_i[U_i])=2\lfloor\frac{n}{3}\rfloor$, so let $c_i: E(D'_i[U_i])\to \{1,2,\ldots,2\lfloor\frac{n}{3}\rfloor\}$ be a proper-coloring of $D'_i[U_i]$.
  Let $D'$ be the (disjoint) union of $D'_1,D'_2,D'_3$.
  We derive $D''$ from $D'$ by lifting each edge of $c^{-1}_i(j)$ to the $j$\textsuperscript{th} vertex of $V_{i+1}\cup V_{i+2}$ (take the indices cyclically). Observe that $D''$ is a simple bipartite graph, whose color classes are still $U$ and $V$, and it is obtained from $D$ via a series of liftings, therefore it is a resolution of $D$.
\end{proof}

\begin{proof}[\normalfont\bfseries Proof of Theorem~\ref{thm:deg}]
  Let us assume that $a\geq b$ and $A=\{v_1,\dots,v_a\}$. By adding edges, if necessary, we may assume that $D$ is semiregular with degrees $\Delta_A$ and $\Delta_B$, where $|E(D)|=a\cdot\Delta_A=b\cdot\Delta_B$. As $D$ is bipartite, by Kőnig's theorem we have $\chi'(D)=\Delta(D)=\Delta_B$, which means that we can split the edges of $D$ into $\Delta_B$ matchings of size $b$, say $M_1,M_2,\ldots,M_{\Delta_B}$.

  \medskip

  We claim that by splitting these matchings appropriately, we can get a partition of the edges of $D$ into matchings $M_1',M_2',\dots,M_a'$, each of size $\Delta_A$. Pick $\Delta_A$ edges of $M_1$ arbitrarily to get $M_1'$ and continue picking sets of $\Delta_A$ edges of $M_1$ that are disjoint from the previously chosen sets, until less than $b/4$ edges of $M_1$ are available. Put the remaining edges into a new $M_i'$; it is easy to see that these edges intersect at most $b/2$ edges of $M_2$, so we can pick some of these edges of $M_2$ to fill up $M_i'$ to the appropriate size. Continue this procedure until less than $b/4$ edges remain in $M_{\Delta_B}$.
  However, as $a=|E(D)|/\Delta_A$, this means that actually all the edges in $M_{\Delta_B}$ are used up as well, thus our claim is proven.

  \medskip

  For each $1\leq i\leq a$, we lift the edges of $M_i'$ to $v_i\in A$. Let us call the resulting demand graph $D'$. In $D'$ there are no multiple edges between $A$ and $B$, $\mu(D'[A])\le 2$, $e_{D'}(v_i,A)\leq2\Delta_A$ and $e_{D'}(v_i,B)=\Delta_A$ for all $v_i\in A$.

  \medskip

  To each edge $e$ with end vertices $\{ v_i,v_j\}\subset A$, we associate a list $L(e)$ of vertices of $B$, to which we can lift $e$ to without creating multiple edges:
  \[L(e)=V(B)\setminus \left(N_{D'}(v_i)\bigcup N_{D'}(v_j)\right).\]
  We have $|L(e)|\geq b-e_{D'}(v_i,B)-e_{D'}(v_j,B) \geq b-2\Delta_A$. By Kahn's theorem (Theorem~\ref{thm:kahn}), $\text{ch}'(D'[A])\leq(1+o(1))\chi'(D'[A])$.
  Furthermore, by Vizing's theorem (Theorem~\ref{thm:vizing}), $\chi'(D'[A])\leq \Delta(D'[A])+\mu(D'[A]) \leq 2\Delta_A+2$. By the assumptions made in the statement of the theorem on $\Delta_A$, we have $\mathrm{ch}'(A)\le |L(e)|$ for each edge $e$ in $E(D'[A])$.
  Thus, there is a proper list edge coloring $c$ which maps each $e\in E(D'[A])$ to an element of $L(e)$. Finally, we lift every edge $e\in E(D'[A])$ to $c(e)$. As we do not create multiple edges between $A$ and $B$, the resulting graph is a resolution of $D$.
\end{proof}

\section{Proof of the edge version (Theorem~\ref{thm:edge})}
We proceed by mathematical induction on $n$. It is easy to check that the result holds for $n=4,5$ by a straightforward case analysis. 

\medskip

Let $A$ and $B$ be the color classes of $D$, each of cardinality $n$. In the induction step we lift some edges in $D$ in such a way that the resulting graph $D'$ is still bipartite with the same color classes and there exists a subset $Z\subset V(D')$ such that
\begin{enumerate}
  \item $|Z\cap A|=|Z\cap B|$ holds,
  \item $\ge |Z|$ edges of $D'$ are incident to vertices of $Z$,
  \item $\Delta(D'[(A\cup B)\setminus Z])\le n-|Z|/2$, and
  \item there are no multiple edges incident to vertices of $Z$ in $D'$.
\end{enumerate}
The first three conditions guarantee that we can invoke the inductive hypothesis on $D'[(A\cup B)\setminus Z]$, to conclude that $D'[(A\cup B)\setminus Z]$ is resolvable. The fourth condition now implies that $D'$ is resolvable as well, which in turn implies the same for $D$.

\medskip

Since we want to keep $D'$ bipartite with the same color classes as $D$, we define the \emph{edge-lifting} of an edge $e\in E(D)$, with end vertices $u\in A$ and $v\in B$, to $xy$, whenever $\{u,v,x,y\}$ are four different vertices and $x\in A$ and $y\in B$: the operation adds a copy of $xy$, $uy$, and $xv$ to $D$ and then deletes $e$. Note that an edge-lifting operation can also be obtained as a composition of two liftings (one to $x$ and then to $y$).

\medskip

Assume now that $n \geq 6$ and let $D$ be a demand graph on $2n-2$ edges. (We may make this assumption on the number of edges by adding edges between two vertices of degree less than $n$ from distinct classes.) Let
\[X=\{v\in A\cup B\ :\ d(v)=n\}.\]
As we have $2n-2$ edges, it is clear that $X$ meets both $A$ and $B$ in at most one vertex, so $|X|\leq 2$. Furthermore, each color class has either at least one isolated vertex or at least two vertices of degree $1$.

\medskip

We distinguish four major cases.

\makeatletter
  \renewcommand{\thesubsection}{\bfseries Case~\arabic{subsection}}
  \renewcommand\subsection{\@startsection{subsection}{2}{\z@}%
  	{9\p@ \@plus 6\p@ \@minus 3\p@}%
  	{3\p@ \@plus 6\p@ \@minus 3\p@}%
  	{\normalfont\normalsize}}

  \renewcommand{\thesubsubsection}{\bfseries Case~\arabic{subsection}.\arabic{subsubsection}}
  \renewcommand\subsubsection{\@startsection{subsubsection}{2}{\z@}%
  	{9\p@ \@plus 6\p@ \@minus 3\p@}%
  	{3\p@ \@plus 6\p@ \@minus 3\p@}%
  	{\normalfont\normalsize}}

  \renewcommand{\theparagraph}{\bfseries Case~\arabic{subsection}.\arabic{subsubsection}.\arabic{paragraph}}
  \renewcommand\paragraph{\@startsection{paragraph}{2}{\z@}%
  	{9\p@ \@plus 6\p@ \@minus 3\p@}%
  	{3\p@ \@plus 6\p@ \@minus 3\p@}%
  	{\normalfont\normalsize}}

\makeatother

\subsection{\texorpdfstring{$u_1,u_2\in A$ and $v_1,v_2\in B$}{u\_1,u\_2∈A and v\_1,v\_2∈B} are four isolated vertices in \texorpdfstring{$A$ and $B$.}{A and B.}}\label{item:twoisolated}

Let $Y=\{v\in A\cup B\ :\ d(v)\ge n-1\}$ and set $Z=\{u_1,u_2,v_1,v_2\}$. Suppose there exists a set $F\subset E(D)$ of four edges, which cover every vertex of $D$ at most twice, cover every element of $Y$ at least once, and cover every element of $X$ exactly twice.
It is easy to see that there is a numbering $F=\{e_1,e_2,e_3,e_4\}$ of these edges, so that edge-lifting $e_1$ to $u_1v_1$, $e_2$ to $u_1v_2$, $e_3$ to $u_2v_2$, and $e_4$ to $u_2v_1$ does not create multiple edges. Therefore, given the existence of $F$, we can invoke the inductive hypothesis and conclude that $D$ is resolvable in $K_{n,n}$.

\medskip

Notice, that
\begin{equation}\label{eq:case1}
  \sum_{v\in Y}d(v)-|E(D[Y])|\le |E(D)|=2n-2,
\end{equation}
and $\Delta(D[Y])\le n$. Depending on the cardinality of $|Y|$, we distinguish 5 subcases.

\subsubsection{$|Y|=4$.} Since $\sum_{v\in Y}d(v)\ge 4(n-1)$, by Equation~(\ref{eq:case1}) we have $2n-2\le |E(D[Y])|$, so actually every edge of $D$ is induced by $Y$. If there is a $C_4$ in $D[Y]$, then the edges of the cycle are a good choice for $F$. Otherwise we can pair the vertices of $Y$ in such a way that the pairs are joined by at least $n-2$ edges each; choose two edges from each pair, and let this set of four edges be $F$.

\subsubsection{$|Y|=3$.} Again, we have $\sum_{v\in Y}d(v)\ge 3(n-1)$ in Equation~\ref{eq:case1}, thus $n-1\le |E(D[Y])|$. Also, $Y$ has exactly one vertex in either $A$ or $B$, therefore $|E(D[Y])|\le \Delta(D[Y])\le n$.
Without loss of generality, we may suppose that $A\cap Y=\{a_1\}$ and $B\cap Y=\{b_1,b_2\}$, and that $e(a_1,b_1)\ge (n-1)/2\ge 2$.
Therefore $e(b_2,V(D)\setminus Y)\ge (n-1)/2\ge 2$ as well.
Choose two edges joining $a_1$ to $b_1$ and two edges joining $b_2$ to $V(D)\setminus Y$, and let this set of four edges be $F$.

\subsubsection{$|Y|=2$.} If both $e(A\cap Y, V(D)\setminus Y)\ge 2$ and $e(B\cap Y, V(D)\setminus Y)\ge 2$, then choose two edges from both sets; this set of four edges is a good choice for $F$. Otherwise $|E(D[Y])|\ge n-2$, therefore there are at most $n+2$ edges incident on $Y$, or in other words, $V(D)\setminus Y$ induces at least $n-4\ge 2$ edges. Choose two edges from both $D[Y]$ and $D[V(D)\setminus Y]$, and let this set of four edges be $F$.

\subsubsection{$|Y|=1$.} There is a vertex $v$ to which $Y$ is joined by at least two edges (there are two isolated vertices in both color classes). The vertex $v$ and $Y$ cover at most $2n-4$ edges, so select two edges not intersecting $v$ and $Y$, plus two edges joining $v$ and $Y$; let this set of four edges be $F$.

\subsubsection{$|Y|=0$.} There are two vertices joined by at least two edges, as otherwise $D$ is the resolution of itself. We can proceed exactly as in the $|Y|=1$ case.

\medskip

\rule{.3333\textwidth}{0.5pt}

\bigskip

From now on, without loss of generality, we may assume that there is at most one isolated vertex in one of the classes.

\subsection{\texorpdfstring{$X$}{X} is empty.}

\subsubsection{We have a vertex $x$ of degree $1$ in one of the classes, say, $A$.} Suppose first, that $y$ is an isolated vertex of $B$: then we may edge-lift an edge $e\in E(D)$, which is not incident to $x$ or to the neighbor of $x$, to $xy$ and let $Z=\{x,y\}$; the four conditions are satisfied. If there are no isolated vertices in $B$, there are at least two degree 1 vertices in it (the sum of the degrees is $2n-2$); let us denote by $y$ one of the two points not joined to $x$. We let $D'=D$ and $Z=\{x,y\}$, and proceed with induction.

\subsubsection{There is no vertex of degree one in $D$.} We must have at least one isolated vertex in each class. Furthermore, the average degree of the remaining vertices in each class is $(2n-2)/(n-1)=2$, so we either have another isolated vertex or every remaining vertex has degree exactly two.

\medskip

Recall that we may assume that there is at most one isolated vertex in one of the classes.

\paragraph{There is a vertex of degree two without multiple edges.} Put this vertex and an isolated vertex from the other class into $Z$, and invoke the inductive argument.

\paragraph{There are two isolated vertices, $a$ and $b$, in one of the classes.} Without loss of generality, we may assume that $a,b\in A$. All but one vertex of $B$ has degree two. We may assume that each of the non-isolated vertices of $B$ have parallel edges, or else the previous case applies. Let $u$ and $v$ be the vertices in $A$ with highest degrees, and let $z$ be a neighbor of $u$ and $w$ be a neighbor of $v$. We edge-lift $uz$ to $aw$, $wv$ to $bz$ and let $Z=\{a,b,z,w\}$.

\paragraph{There is exactly one isolated vertex in each of $A$ and $B$.} We may assume that every remaining vertex has degree two and is the endpoint of two parallel edges. In this case, let $a_1,\dots,a_n$ and $b_1,\dots,b_n$ be the vertices of $D$, with $a_i$ and $b_i$ connected by two edges for each $1\leq i\leq n-1$ and $a_n$ and $b_n$ isolated. In this setting, we construct a resolution of $D$ by edge-lifting a copy of the edge $a_i b_i$ to $a_{i+1}b_{i+2}$ for each $1\leq i\leq n-2$, and edge-lifting a copy of $a_{n-1}b_{n-1}$  to $a_{n} b_1$.

\subsection{\texorpdfstring{$|X|=1$.}{|X|=1.}}

Let $z\in A$ be the only vertex of degree $n$ in $D$. Notice that, in this case, there is no vertex of degree $n-1$ in $A$ and there exists at least one isolated vertex in $A$. Let us call it $v$.

                \subsubsection{There is a vertex $u$ of degree $1$ in $B$.}
                We have two cases: if it is joined to $z$, we edge-lift a copy of an edge $xy$ independent from $uz$ to $uv$ and let $Z=\{u,v\}$. If $u$ is not joined to $z$, we simply edge-lift any edge incident on $z$ to $uv$ and let $Z=\{u,v\}$.

                \subsubsection{There is no vertex of degree $1$ in $B$.}
                There must be an isolated vertex $u$ in this class, and the average degree of the remaining vertices is $(2n-2)/(n-1)=2$. Therefore, either every remaining vertex has degree exactly two or there is another isolated vertex in $B$.

                        \paragraph{Every vertex in $B$ except $u$ has degree two.} Either one of them has no adjacent multiple edges or the neighborhood of each of them consists of two parallel edges. In the first case, let $x$ be a vertex without multiplicity. We simply edge-lift an edge of $z$ to $uv$ and let $Z=\{x,v\}$. In the latter case, the degree of each vertex is even, and, as the average degree of the vertices in $A/\{v,z\}$ is $1$, we must have another isolated vertex $v'$ in $A$. Let $a$ be a neighbor of $z$ and $b$ be a vertex of $B$ not joined to $z$, let $z'$ be its neighbor. We edge-lift a copy of $az$ to $bv$, $bz'$ to $av'$ and let $Z=\{v,v',a,b\}$.

                        \paragraph{There is another isolated vertex $u'$ in $B$.} The remaining vertices of $A$ have average degree $(n-2)/(n-2)=1$, so all of them have degree one (recall that there is at most one isolated vertex in one of the classes). In the first case, just take a neighbor $x$ of $z$ that has a non-neighbor $y$ of degree one in $A$ (it does exist because $z$ has at least two neighbors). Edge-lift the edge $zx$ to $uy$ and let $Z=\{u,y\}$.

\subsection{\texorpdfstring{$|X|=2$.}{|X|=2.}}
Let $z_1\in A$ and $z_2\in B$ be the vertices of degree $n$. Notice that $z_1$ and $z_2$ must be joined by at least two edges and that there is no other vertex of degree $n$ or $n-1$ in $D$. Furthermore, in each class, we must have an isolated vertex, $v_1\in A$ and $v_2\in B$, and the average degree of the remaining vertices is $(n-2)/(n-2)=1$, so in each class either we have another isolated vertex or all the remaining vertices have degree one.

\medskip

Recall that there is at most one isolated vertex in one of the classes, say $B$. All vertices except $z_2$ and $v_2$ have degree one, then either we have a vertex $x$ of degree one which is not joined to $z_1$, or $z_1$ is joined to every vertex of positive degree in $B$. In the first case, edge-lift a copy of $z_1z_2$ to $v_1x$ and let $Z=\{x,v_1\}$. In the latter case, the neighborhood of $z_1$ consists of $n-2$ simple edges connecting it to the vertices of degree one in $B$ and one double edge joining $z_1$ and $z_2$.
Simply edge-lift one copy of this double edge to $v_1v_2$ and let $Z=\{z_1,v_2\}$ (as $z_1$ has no multiple edges now).

\bigskip

Our case analysis is now complete, as is the proof of Theorem~\ref{thm:edge}.

\phantomsection
\bibliography{ref}

\end{document}